\documentclass[11pt,reqno]{amsart}
\usepackage{amsmath,amsfonts,amssymb,amsthm,amscd}

\makeindex

\usepackage{amssymb}
\usepackage[latin1]{inputenc}
\usepackage{ epsfig,epsf}
\usepackage{enumerate}
\usepackage{indentfirst}
\usepackage{euscript}
\usepackage{graphicx}
\usepackage{amsmath}
\usepackage{amsfonts}
\usepackage{amssymb}
\makeindex

\makeindex
 \theoremstyle{plain}
\newtheorem{proposition}{Proposition}
\newtheorem{lemma}{Lemma}

\newtheorem{remark}{Remark}
\theoremstyle{definition}

\newtheorem{problem}{Problem}

 \title[   Lines of Curvature near Singular End Points]{  Lines of Principal Curvature near Singular
  End Points of Surfaces in ${\mathbb R}^3$}

\author[J. Sotomayor]{Jorge Sotomayor}
\author[R. Garcia]{Ronaldo Garcia}

\address{Instituto de Matem\'{a}tica e Estat\'{\i}stica \\Universidade de S\~{a}o
Paulo
\\Rua do Mat\~{a}o 1010,  Cidade Universit\'{a}ria  \\CEP 05508-090, S\~{a}o Paulo, S.P., Brazil
}

\address{Instituto de Matem\'{a}tica e Estat\'{\i}stica \\
Universidade Federal de Goi\'as \\CEP 74001-970, Caixa Postal 131
\\Goi\^ania, GO, Brazil}

 \thanks{2000 {\it Mathematics Subject Classification}. Primary 53A05, 34C23; Secondary 58K25.\\
{\it Key words and phrases}. Principal curvature lines, inflexion
singular  end points.\\
This work was done under the projects PRONEX/CNPq/MCT - grant
number
 66.2249/1997-6,
 and
 CNPq/PADCT - grant number 620029/2004-8.  It was also
 partially supported by CNPq Grant  476886/2001-5 and FUNAPE/UFG}

\begin{document}

  \begin{abstract}
In this paper are studied the nets of principal curvature lines on
surfaces embedded in Euclidean $3-$space near their end  points,
at which the surfaces tend to infinity.

This is a natural complement and extension to smooth surfaces  of
the work of Garcia and Sotomayor (1996), devoted to the study of
principal  curvature nets  which are structurally stable --do not
change topologically-- under small perturbations on the
coefficients of the equations defining  algebraic surfaces.

 This paper goes one step further and classifies the patterns of the most common and stable  behaviors at the ends,
present also  in generic  families of   surfaces depending on  one-parameter.
\end{abstract}
\maketitle

\section{ Introduction}\label{sec:1}

A  {\it   surface} of {\it smoothness class} $C^k $ in Euclidean $
(x, y, z)$-space $\mathbb R ^3$  is defined by the  variety
$A(\alpha)$ of zeros of a  real function $\alpha $ of class $C^ k$
in  $\mathbb R^3$. The exponent $k$  ranges among the positive
integers as well as on  the symbols  $\infty$, $\omega$ (for
analytic) and   $a(n)$ (for algebraic of degree $n$).

In the  class    $C^{a(n)}$  of algebraic surfaces of degree $n$,
we have $\alpha= \sum \alpha_h , \,h=0,\,1,\,2,...,n\,$, where
$\alpha_h$   is  a   homogeneous   polynomial  of  degree $h\,\,$
with real coefficients:  $\alpha_h = \sum a_{ijk}x^i y^j z^k,
\,\,i+j+k=h.$

The space $\mathbb R^3$  will be endowed with the Euclidean metric
$ ds^2 =  dx^2\,  + dy^2 \, + dz^2 $ also denoted  by $ <  ,  > $,
and  with  the  positive orientation induced by the volume form
$\Omega = dx\wedge dy\wedge dz$.

An  {\it end  point}  or {\it point  at infinity} of $A(\alpha) $
is a point in  the unit sphere $\mathbb S^2$, which is the limit
of  a sequence  of  the  form $p_n /|p_n|$, for $ p_n$  tending to
infinity in $A(\alpha). $

The  {\it end  locus}, $E(\alpha)$,
 of $A(\alpha)$ is
the collection  of its end  points. This
set is a geometric measure of the non-compactness of the surface
and describes how it tends to infinity.

A surface $A(\alpha)$ is said to be {\it  regular
 at}  $p\in E(\alpha ) $
if in a neighborhood of $p$, $ E(\alpha ) $ is a
  regular smooth  curve
in  $\mathbb S^2.$  Otherwise, $p$  is said to be a {\it
critical end point} of  $A(\alpha)$.

For the class $a(n)$,  $E(\alpha) $ is contained in  the
algebraic  curve $E_n(\alpha) = \{p \in \mathbb S^2 ; \; \alpha_n (p)=
0\}.$  The regularity of $E(\alpha)$
is equivalent to that of $E_n(\alpha)$.

The  gradient  vector field of $\alpha,$
will be denoted by  $\nabla\alpha  =  \alpha_x \partial / \partial
x  + \alpha_y \partial / \partial y  + \alpha_z \partial /
\partial z ,$ where  $\alpha_x =\partial\alpha / \partial x   $, etc.

The   zeros  of  this  vector field are called {\it critical
points} of $\alpha$; they  determine  the  set $C(\alpha).$ The
regular part of $A(\alpha) $ is  the  smooth  surface $ S(\alpha)
= A(\alpha)\setminus C(\alpha).$  When $C(\alpha)$ is disjoint
from $A(\alpha)$, the surface $ S(\alpha) = A(\alpha) $ is called
{\it regular}. The {\it orientation}  on  $S(\alpha) $ will  be defined by
taking the gradient $\nabla\alpha $ to  be  the {\it positive
normal}. Thus  $ A(-\alpha)$   defines   the  same
surface  as $A(\alpha)$  but endowed   with   the   opposite
orientation on $S(-\alpha).$

 The  Gaussian  normal  map  $N$,  of
$S(\alpha)$ into  the sphere $\mathbb  S^2$,  is defined by the
unit vector in the direction  of  the  gradient:
$N_\alpha=\nabla\alpha/|\nabla\alpha|.$ The eigenvalues
$-k_\alpha^1  (p)$ and  $-k_\alpha^2  (p)$ of   the   operator
$DN_\alpha (p),$  restricted to $T_p S(\alpha),$ the tangent space
to the surface  at $ p$, define the principal  curvatures,
$k_\alpha^1  (p)$ and  $k_\alpha^2 (p)$   of  the surface at the
point $p.$ It will be  always  assumed  that  $ k_\alpha^1  (p)
\leq   k_\alpha^2 (p).$

The  points  on $ S(\alpha)$ at which the principal curvatures
coincide, define the set $U(\alpha)$ of  {\it umbilic  points} of
the surface   $A(\alpha). $   On $S(\alpha)\setminus U(\alpha)$,
the eigenspaces of $DN_\alpha$, associated   to $-k_\alpha^1 $ and
$-k_\alpha^2 $ define line fields $L_1 (\alpha) $ and $
L_2(\alpha)$, mutually orthogonal,  called respectively {\it
minimal} and  {\it maximal principal  line fields}  of the surface
$A(\alpha).$ The  smoothness class of these line fields is
$C^{k-2}$, where $k-2 = k$
 for  $k= \infty, \; \omega$ and  $a(n)-2 = \omega $.

The  maximal  integral  curves  of  the  line  fields $ L_1
(\alpha) $ and  $L_2 (\alpha) $ are called respectively  the lines
of  minimal  and  maximal  principal curvature,  or simply the
{\it principal  lines} of $A(\alpha)$.

 What was  said above
concerning the definition of these lines is equivalent   to
require   that  they  are non  trivial  solutions of   Rodrigues'
differential equations:

\begin{equation}\label{eq:11}
DN_\alpha (p)dp + k_\alpha^i  (p)dp = 0,\;\;\;\;\; <N (p),dp>
=0,\;\; i=1,2.
 \end{equation}

\noindent where $p=(x ,y ,z ),\;  \alpha(p)=0,$ $ dp= dx
\partial /\partial x \; + dy
\partial /\partial y \; + dz
\partial /\partial z.$
See \cite{Sp, St}.

After  elimination  of $ k^i_\alpha ,\; i=1,2$, the first two  equations  in
(\ref{eq:11}) can be written as the following single implicit
quadratic equation:

\begin{equation}\label{eq:11l}
 <DN_\alpha (p)dp\wedge N (p),dp> = [DN_\alpha (p)dp,N_\alpha (p),dp]=0. \end{equation}

The
left (and mid term)    member of this equation is the
{\it geodesic torsion} in the direction of $dp$. In terms of a
local parametrization $\overline{\alpha}$ introducing coordinates
$(u,v)$ on the surface, the equation of lines of curvature in
terms of the coefficients  $(E,F,G)$ of the first and  $(e,f,g)$
of the second fundamental forms is, see \cite{Sp, St},

\begin{equation}\label{eq:11p}
[Fg-Gf]dv^2 + [Eg-Ge]dudv + [Ef-Fe]du^2 = 0.
 \end{equation}

The  net  $F(\alpha)=(F_1(\alpha),F_2 (\alpha)) $ of  orthogonal
curves  on $ S(\alpha)\setminus U(\alpha),$ defined by the
integral  foliations  $F_1 (\alpha)$  and  $ F_2 (\alpha)$  of the
line fields $ L_1 (\alpha)$ and $L_2 (\alpha)$,  will  be called
{\it the principal  net }  on $A(\alpha)$.

The study of families of principal  curves and their umbilic
singularities
on   immersed surfaces
was initiated
by
Euler, Monge, Dupin and
Darboux,
 to mention only a few.
See \cite{Da, Mo}  and   \cite{gs1, Sp,  St} for references.

Recently this classic subject acquired new vigor   by the introduction
of ideas coming from Dynamical Systems and the  Qualitative Theory
of Differential Equations.
See the works  \cite{gs1}, \cite {gas1}, \cite{gs4}, \cite{G-S4}       of
 Gutierrez, Garcia and Sotomayor
on the structural stability, bifurcations and genericity of
principal curvature lines and their umbilic and critical
singularities  on compact surfaces.

The scope of the subject was broadened by the extension of the works  on structural  stability to
other  families of curves
 of classical geometry. See \cite {ggs},
for the  asymptotic lines and  \cite {amean, tou, sbm, gm} respectively  for  the
arithmetic, geometric,
 harmonic and general   mean curvature lines. Other pertinent directions of research
 involving implicit differential equations arise from Control and Singularity Theories,
 see Davydov \cite{dav} and Davydov, Ishikawa, Izumiya and Sun \cite{dav2}.

In \cite{gasalg}  the authors studied the behavior of the lines of
curvature on algebraic surfaces, i.e.  those of $C^{a(n)}$, focusing particularly their generic and stable
patterns  at  end points. Essential for this study was  the
operation of compactification of algebraic surfaces and their equations (\ref{eq:11l}) and
 (\ref{eq:11p}) in $\mathbb R^3$  to obtain compact ones in $\mathbb S^3$. This step is  reminiscent of the Poincar\'e
compactification of polynomial differential equations \cite {ps}.

In this paper the  study in \cite{gasalg}  will be extended to the broader and more flexible case of $C^{k}$-smooth surfaces.

As mentioned above, in the case of algebraic surfaces studied in   \cite{gasalg}, the ends
are the algebraic curves defined by the zeros, in the Equatorial Sphere $\mathbb S^2$ of  $\mathbb S^3$, of the highest degree homogeneous part $\alpha _n$ of the polynomial  $\alpha $. Here, to make the study of the principal nets at  ends of smooth surfaces tractable by methods of Differential
Analysis,
we follow an inverse procedure, going from compact
 smooth surfaces
in $\mathbb S^3$
to surfaces
in $\mathbb R^3$.
This restriction on the class of surfaces studied in this paper  is explained in Subsection \ref{sec:pre}.

The new results of this paper on the  patterns of principal nets
at end points  are  established in Sections \ref{sec:end1} and
\ref{sec:end2}. Their meaning for the Structural Stability and
Bifurcation Theories of Principal Nets is discussed in Section
\ref{sec:fin}, where a pertinent problem is proposed. The essay
\cite{hc} presents  a historic overview of the subject and reviews
other problems left open.

\subsection {Preliminaries} \label{sec:pre}

Consider the  space ${\mathcal A}_c ^k$  of  real valued functions
$\alpha^c $ which are  $C^k$-smooth
 in   the   three
dimensional   sphere $\mathbb S^3 =\{|p|^2 +|w|^2 =1\}$ in $
\mathbb R ^4$, with coordinates $p=(x, y ,z)$ and $w$. The meaning
of the exponent  $k$ is the same as above and $a(n)$ means
polynomials of degree $n$ in four variables of the form $\alpha^c
= \sum \alpha_h w^h , \,h=0,\,1,\,2,...,n$, with $\alpha_h$
homogeneous of degree $h$ in $(x,y,z)$.

The {\it
equatorial} sphere in  $\mathbb S^3$  will be
 $\mathbb S^2 =\{(p,w): \; |p|=1, \; w=0\}$ in $
\mathbb R^3$. It will be endowed with the positive orientation defined  by the outward normal.
 The  northern
hemisphere of  $ \mathbb S^3$ is defined by ${\mathbb H}^+
=\{(p,w)\in \mathbb S^3: \; w>0\}$.

The surfaces $A(\alpha)$ considered in this work  will   defined
in terms of functions $\alpha ^c    \in {\mathcal A}_c ^k$ as
$\alpha= \alpha ^c  {\circ}  \mathbb P$,  where $\mathbb P$ is the
{\it central projection} of  $\mathbb R^3$,  identified with the
tangent plane at the north pole $\mathbb T_w^3$, onto $\mathbb H
^{+}$, defined by:

$$\mathbb P (p)   = (p/(|p|^2 +1)^{1/2}   ,\;\; 1/(|p|^2 +1)^{1/2}).$$

For future reference, denote by $\mathbb T_y^3$ the tangent plane to $\mathbb S^3$ at the point $(0,1,0,0)$, identified with $\mathbb R^3$ with orthonormal coordinates $(u,v,w)$, with $w$ along the vector $\omega=(0,0,0,1)$. The central projection $\mathbb Q$ of $\mathbb T_y^3$ to $\mathbb S^3$ is such that ${\mathbb P}^{-1}\circ \mathbb Q:\mathbb T_y^3\to \mathbb T_w^3$ has the coordinate expression $(u,v,w)\to (u/w, v/w,1/w)$.

 For $m \leq k$, the  following expression defines uniquely, the functions involved:

\begin{equation} \label{der}
\alpha ^c  (p,w)=  \sum w ^{j} \alpha_j^{c} (p)\; + \; o(|w|^m ),
\,j=0,\,1,\,2,...,m.
\end {equation}

In the algebraic case ($k=a(n)$) studied in \cite {gasalg},
$\alpha ^c  =  \sum w ^{n-h} \alpha_h , \,h=0,\,1,\,2,..,n$, where
the obvious correspondence  $ \alpha_h = \alpha_{n-h} ^{c} $
holds.

The end points of $A(\alpha)$, $E(\alpha)$,  are contained in
$E(\alpha ^c ) =\{\alpha_0^{c} (p)=0\}.$

At a  {\it regular  end  point}  $p$ of $ E(\alpha)$  it
 will be  required that $\alpha_0^{c}$ has a regular zero,  i.e. one with non-vanishing derivative i.e.
 $\nabla \alpha_0^{c} (p) \neq 0.$
At regular end points, the end locus is oriented by the positive
unit normal  $\nu(\alpha) = \nabla\alpha_0^{c}
/|\nabla\alpha_0^{c} |. $
This defines the positive unit tangent
vector along $ E(\alpha)$, given at $p$ by $\tau(\alpha) (p) =
p\wedge \nu(\alpha) (p). $ An end point is called {\it critical}
if it is not regular.

A  regular point $p\in E(\alpha)$  is    called  a {\it biregular  end point} of $
A(\alpha)$ if  the geodesic curvature, $ k_g$, of
the  curve $E(\alpha) $  at $p$, considered  as  a  spherical
curve,  is different  from zero; it is called
an {\it inflexion end point} if $k_g$  is equal to zero.

  When the surface  $A(\alpha)$ is regular at infinity, clearly $E(\alpha) = E(\alpha ^c) =\{ \alpha_0^{c} (p)=0\}.$

The analysis in sections \ref{sec:end1} and \ref{sec:end2}   will  prove
that there is a natural   {\it extension}
  $F_c (\alpha)=(F_{c1}(\alpha),F_{c2} (\alpha)) $   of  the net
  ($\mathbb P (F_1(\alpha)), \mathbb P (F_2(\alpha))$  to $A_c (\alpha) = \{\alpha ^c =0\}$,  as  a net of class $C^{k-2}$, whose  singularities in $E(\alpha)$ are  located at
    the {\it inflexion}
      and
{\it critical} end points of $A(\alpha)$. This is done by means of
special charts used to extend the quadratic differential equations
that define ($\mathbb P (F_1(\alpha)), \mathbb P (F_2(\alpha))$ to
a full neighborhood in $A_c (\alpha)$ of the arcs of biregular
ends. The differential equations are then  extended to a full
neighborhood of  the singularities. See Lemma \ref{lem:d1}, for
regular ends, and Lemma \ref{lem:d2}, for critical ends.

 The main contribution of this paper consists  in
  the resolution  of singularities of the extended differential equations, under suitable genericity hypotheses on $\alpha^c$.
This is done in sections \ref{sec:end1} and \ref{sec:end2}.  It  leads
to
eight  patterns of principal nets  at  end points. Two of them -- elliptic and hyperbolic inflexions-- have  also  been  studied in the case of  algebraic surfaces \cite {gasalg}.

\section{Principal Nets at Regular End Points }\label{sec:end1}

\begin{lemma} \label{lem:c1}
Let $p$ be a regular end point of $A(\alpha)$,
$\alpha=\alpha^c\circ \mathbb P$. Then there is   a mapping
$\overline{\alpha}$ of the form

$\overline{\alpha}(u,w)=(x(u,w), y(u,w), z(u,w))$, $w >0$,  defined by

\begin{equation} \label{eq:pro}
x(u,w)= \frac{u}{w},\;\;\;\;\;\;\;  y(u,w)=  \frac{h(u,w)}w,\;\;\;\;\;\;\;
z(u,w)= \frac{1}{w}.
\end{equation}
which parametrizes the surface $A(\alpha)$ near $p$, with

\begin{equation} \label{eq:rep}
\aligned
 h(u,w)=&
 k_0 w+\frac 12 au^2+ b  u w+\frac 12  c w^2\\
 +&\frac 16(a_{30} u^3+3a_{21}u^2w+3 a_{12}uw^2+a_{03}w^3)\\
 +& \frac 1{24}( a_{40}u^4+4
 a_{31}u^3w+6a_{22}u^2w^2\\
+&4a_{13}uw^3+a_{04}w^4)+ h.o.t
\endaligned
\end{equation}

\end{lemma}

\begin{proof}
With no lost of generality,  assume  that the
regular end  point $p$ is located at $(0,1,0,0)$,
the unit tangent vector to the regular end curve is $\tau =(1,0,0,0)$
and the positive normal vector is $\nu = (0,0,1,0)$.
Take orthonormal coordinates $u,v,w$ along $\tau , \nu, \omega =(0,0,0,1)$
on the tangent space, ${\mathbb T }_y^3$ to $\mathbb S ^3$ at $p$.
Then the composition ${\mathbb P}^{-1}\circ \mathbb Q$  writes as
 $x=u/w, \; y=v/w, \; z=1/w$.

Clearly the surface $A(\alpha^c)$ near  $p$  can be parametrized
by the
 central projection into $\mathbb S ^3$ of the graph of a $C^k$ function of the form
$v=h(u,w)$ in ${\mathbb T}_y^3$, with $h(0,0)=0$ and $h_u(0,0)=0$.
This means that the surface $A(\alpha)$, with $\alpha =
\alpha^c\circ \mathbb P ^{-1}$ can be parametrized in the form
(\ref{eq:pro}) with $h$ as in (\ref{eq:rep}).

\end{proof}

\begin{lemma} \label{lem:d1} The differential equation (\ref{eq:11p}) in the chart $\overline{\alpha}$ of Lemma \ref{lem:c1}, multiplied by $w^5\sqrt{EG-F^2}$, extends to a full domain of the chart $(u,w)$  to one   given by
\begin{equation} \label{eq:lmn} \aligned
   L&dw^2+Mdudw+Ndu^2=0, \\
 L  =&-b-a_{21}u-a_{12}w -(c+a_{22})uw\\
-&(b+\frac 12a_{31}) u^2-\frac 12
 a_{13}w^2 + h.o.t.\\
M =& -a-a_{30}u-a_{21}w-\frac 12(2a+a_{40}) u^2\\
-&a_{31}uw+\frac 12(2c-a_{22})w^2+h.o.t.\\
N  =& w[a u  +b w +a_{30} u^2 +2a_{21}uw +a_{12}w^2 +h.o.t.]
\endaligned
\end{equation}
\noindent where the coefficients are of   class $C^{k-2}$.
\end{lemma}

\begin{proof}
The coefficients of  first fundamental form of $\overline{\alpha}$
in (\ref{eq:pro}) and (\ref{eq:rep})  are given by:

$$\aligned
 E(u,w) =& \frac{ 1+h_u^2 }{w^2}\\
F(u,w)=& \frac{ h_u(  w  h_w -h) -u}{w^3} \\
G(u,w)=& \frac{ 1+ u^2+  (w h_w - h )^2  }{w^4}
\endaligned
$$
 The coefficients of the second fundamental
form of $\overline{\alpha}$ are:

$$\aligned
 e(u,w) =&  \frac{ h_{uu }}{w^4\sqrt{EG-F^2}}, \\
f(u,w)=& \frac{ h_{uw}  }{w^4\sqrt{EG-F^2}}, \\
g(u,w)=& \frac{h_{ww}}{w^4\sqrt{EG-F^2}}
\endaligned
$$
\noindent where $e=[\overline{\alpha}_{uu}, \overline{\alpha}_u,
\overline{\alpha}_w]/|\overline{\alpha}_u\wedge
\overline{\alpha}_w|, $  $ f=[\overline{\alpha}_{uw},
\overline{\alpha}_u,
\overline{\alpha}_w]/|\overline{\alpha}_u\wedge
\overline{\alpha}_w|  $ and $g=[\overline{\alpha}_{ww},
\overline{\alpha}_u,
\overline{\alpha}_w]/|\overline{\alpha}_u\wedge
\overline{\alpha}_w|. $

 The  differential equation of curvature lines (\ref{eq:11p}) is given by
 $Ldw^2+Mdudw+Ndu^2=0, $
where $L=Fg-Gf, \;\; M=Eg-Ge $ and $ N=Ef-Fe$.

These coefficients, after   multiplication by $w^5\sqrt{EG-F^2}$,
keeping the same notation,  give  the expressions in
(\ref{eq:lmn}).
\end{proof}

The differential equation (\ref{eq:lmn}) is non-singular, i.e.,
defines  a regular net of transversal curves if $a \neq 0$.  This
will be seen in item a) of next proposition. Calculation expresses
$a$ as   a non-trivial factor of  $k_g$.

The singularities  of equation (\ref{eq:lmn}) arise when $a=0$;
they will be resolved in item b), under the genericity hypothesis
$a_{30}b\ne 0$.

\begin{proposition}\label{prop:he} Let  $\overline{\alpha} $ be
as in Lemma \ref{lem:c1}. Then the end locus is  parametrized by
the regular curve $v=h(u,w),\; w=0$.

\begin{itemize}
 \item[a)]  At a
 biregular  end  point,  i.e.,  regular and non inflexion, $ a\ne 0, $ the
  principal net is as
 illustrated in Fig. \ref{fig:he}, left.

\item[b)]  If  $p  $ is  an  inflexion,  bitransversal  end point,
i.e., $ \beta(p)= a_{30} b\ne 0$, the principal net is as
illustrated  in  Fig.  \ref{fig:he},    hyperbolic $\beta <0$,
center,   and elliptic $\beta  >0$, right.

\begin{figure}[htbp]

\par
\begin{center}
 \includegraphics[angle=0, width=10.5cm]{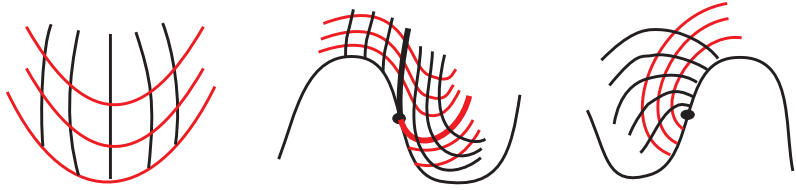}
\end{center}
\caption{ \label{fig:he} Curvature lines near  regular end
 points: biregular (left ) and  inflexions:  hyperbolic (center)  and elliptic (right). }
\end{figure}
\end{itemize}
\end{proposition}

\begin{proof}
Consider the implicit differential equation
\begin{equation}\label{eq:kgn0}\aligned
{\mathcal F}(u,w,p) =& -(  b+a_{21}u+a_{12}w +h.o.t. )p^2\\
 -& (a+a_{30}u+a_{21}w+ h.o.t.)p \\
+&   w( a u + b w+h.o.t. ) =0. \endaligned \end{equation}

 The Lie-Cartan
line field  tangent to  the surface  ${\mathcal F}^{-1}(0)$ is defined by
 $X=({\mathcal F}_p, p {\mathcal F}_p, -({\mathcal F}_u+p{\mathcal
F}_w)) $ in the chart $p=dw/du$ and $Y=(q{\mathcal F}_q, {\mathcal
F}_q, -(q{\mathcal F}_u+ {\mathcal F}_w))$ in the chart $q=du/dw$.
Recall that the integral curves of this line field projects to the
solutions of the implicit differential equation (\ref{eq:kgn0}).

If $a\ne 0$, ${\mathcal F}^{-1}(0)$ is a regular surface,
$X(0)=(-a,0,0) \ne 0$ and $Y(0)=(b,a,0)$.  So by the Flow Box theorem the two
principal foliations are regular and transversal near $0$. This ends the proof of item a).

If $a  =0$, ${\mathcal F}^{-1}(0)$ is a quadratic cone  and
$X(0)=0$. Direct calculation shows that
$$\aligned DX(0)=& \left( \begin{matrix}  -a_{30} &- a_{21} &-2b \\
0&0&0\\ 0&0& a_{30}\end{matrix} \right )\endaligned $$

Therefore $0$ is a saddle point with non zero eigenvalues
$-a_{30}$ and $a_{30}$ and the associated eigenvectors are
$e_1=(1,0,0)$ and $e_2=(b, 0, -a_{30})$.

 The  saddle separatrix  tangent to $e_1$ is parametrized by $w=0$
 and has the following parametrization $(s, 0, 0)$.
The   saddle separatrix  tangent to $e_2$ has the following
parametrization:

$$u(s)=s+ O(s^3), \; w(s) = - \dfrac{a_{30}}b \frac{s^2}2 + O(s^3), \;  p(s)= - \dfrac{a_{30}}b s +  O(s^3).$$

If  $a_{30} b<0 $ the projection $(u(s), w(s))$ is contained in
the semiplane $w\geq 0$.  As the  saddle separatrix is
transversal to the plane $\{p=0\}$ the phase portrait of $X$ is as
shown in the Fig. \ref{fig:cone} below. The projections of the
integral curves in the plane $(u,w)$ shows the configurations of the
principal lines near the inflexion point.

\begin{figure}[htbp]
\par
\begin{center}
\includegraphics[angle=0, width=10.5cm]{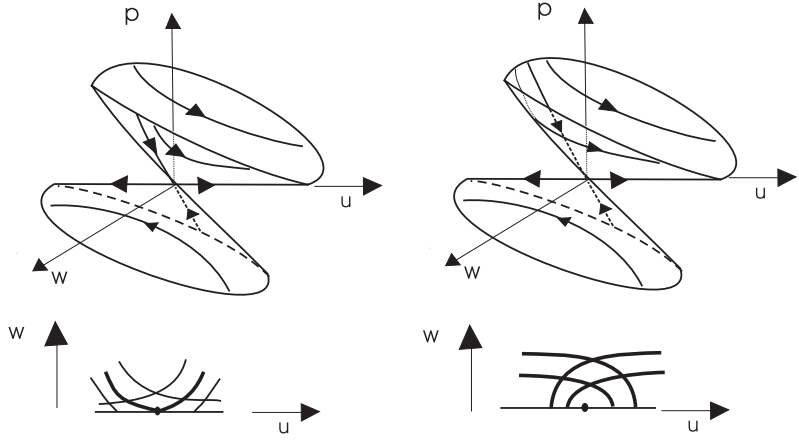}
\end{center}
\caption{ \label{fig:cone} Phase portrait of $X$ near the singular
point of saddle type }
\end{figure}
\end{proof}

\begin{figure}[htbp]

\par
\begin{center}
\includegraphics[angle=0, width=6.5cm]{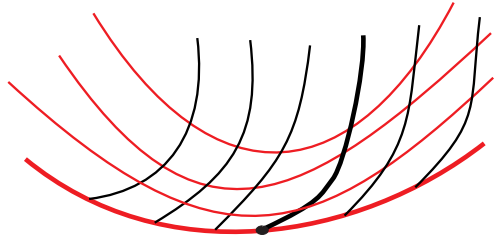}
\end{center}
\caption{ \label{fig:kge} Curvature lines near a
hyperbolic-elliptic inflexion end
point }
\end{figure}

\begin{proposition}\label{prop:kge} Let  $\overline{\alpha} $ be
as in Lemma \ref{lem:c1}. Suppose that, contrary to the hypothesis
of Proposition \ref{prop:he},   $a=0, \; a_{30}=0,$ but
$ba_{40}\ne 0$ holds.

The  differential equation  (\ref{eq:lmn}) of   the principal
lines in this case has the coefficients given by:

 \begin{equation} \label{eq:kgo} \aligned L(u,w)=& -[b+
a_{21}u+a_{12}w+\frac 12 ( 2 b+a_{31}) u^2\\
+&(c+a_{22}) u w+\frac
12a_{13} w^2
 +h.o.t.]
  \\
   M(u,w)=& -[  a_{21}w+\frac 12 a_{40}u^2\\
+&a_{31}uw+\frac 12( a_{22}-2c)w^2
  +h.o.t.]\\
  N(u,w)= &  w^2(b+2a_{21}u +a_{12}w +h.o.t.)\endaligned
  \end{equation}

The principal net is as illustrated in Fig. \ref{fig:kge}.

\end{proposition}

\begin{proof} From equation (\ref{eq:lmn}) it follows the expression
of equation (\ref{eq:kgo}) is  as stated.
  In a neighborhood of $0$ this differential equation factors in  to the product of
  two
  differential
  forms
  $ X_+(u,w)= A(u,w)dv- B_+(u,w)du $  and $ X_-(u,w)= A(u,w)dv- B_-(u,w)du $,
  where $A(u,w)=2L(u,w)$
  and
  $B_{\pm}(u,w)= M(u,w)\pm \sqrt{(M^2-4LN)(u,w)}$.
  The function $A$ is of class $C^{k-2}$ and
  the functions $B_\pm $ are Lipschitz. Assuming $a_{40}>0$, it follows that
$A(0)=-2b\ne 0$,  $B_-(u,0)=0$ and $B_+(u,0)=a_{40}u^2+h.o.t.$ In the case $a_{40}<0$  the analysis is similar,
exchanging $B_-$ with $B_+$.

Therefore, outside the point $0$, the integral leaves of $X_+$ and $X_-$ are transversal.
Further calculation shows that the integral curve  of $X_+$ which
pass through $0$ is parametrized by $(u, -\frac{a_{40}}{6b} u^3 +h.o.t.)$.

This shows that the principal foliations are extended to regular foliations which however fail to be a net a single point of cubic contact.
This is illustrated in Fig. \ref{fig:kge} in the case $a_{40}/b <0$. The case $a_{40}/b >0$ is the mirror image of Fig. \ref{fig:kge}.
\end{proof}

\begin{proposition}\label{prop:kgt}
Let  $\overline{\alpha} $ be as in Lemma \ref{lem:c1}. Suppose
that, contrary to the hypothesis of Proposition \ref{prop:he},
$a=0, \; b=0,$ but $ a_{30} \ne 0$ holds.

The  differential equation   of   the principal lines in this
chart is given by: \begin{equation} \label{eq:kgt} \aligned -&[
a_{21}u+a_{12}w+\frac 12 a_{31} u^2+(c+a_{22}) u w+\frac 12a_{13}
w^2
 +h.o.t.]
  dw^2\\
   -&[ a_{30}u+ a_{21}w+\frac 12 a_{40}u^2+a_{31}uw+\frac 12( a_{22}-2c)w^2
  +h.o.t.]dudw\\
  + & w [ ( a_{30}u^2  +  2a_{21}uw +a_{12}w^2)+h.o.t. ]du^2=0.\endaligned
  \end{equation}

\begin{itemize}
 \item[a)]  If  $(a_{21}^2-a_{12}a_{30})<0 $  the  principal net is as
 illustrated in Fig. \ref{fig:kgt} (left).

\item[b)]  If  $(a_{21}^2-a_{12}a_{30}) > 0 $  the principal net
is as illustrated in Fig. \ref{fig:kgt}(right).

\end{itemize}

\begin{figure}[htbp]

\par
\begin{center}
\includegraphics[angle=0, width=10.5cm]{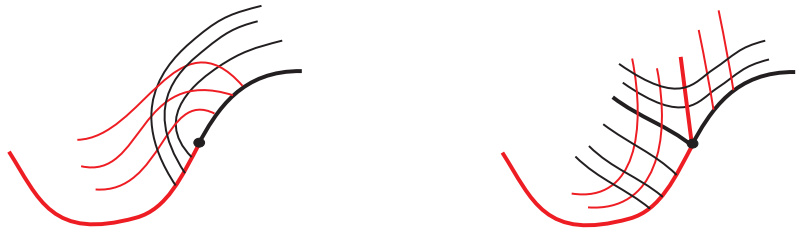}
\end{center}
\caption{ \label{fig:kgt} Curvature lines near an
umbilic-inflexion end point }
\end{figure}

\end{proposition}

\begin{proof}
Consider the Lie-Cartan line field defined by   $$
X=  ({\mathcal F}_p,p {\mathcal F}_p, -({\mathcal
F}_u+p{\mathcal F}_w))$$
 on the singular surface ${\mathcal
F}^{-1}(0)$, where
$$ \aligned  {\mathcal F}(u,w,p)=&- [ a_{21}u+a_{12}w+
  h.o.t.]p^2- [ a_{30}u+ a_{21}w +h.o.t.]p\\
  +  &w[  ( a_{30}u^2  +  2a_{21}uw +a_{12}w^2)+h.o.t. ] =0.\endaligned $$

The singularities of $X$ along the projective line (axis $p$) are
given by the polynomial equation $ p (a_{30}+2 a_{21}p+a_{12}
p^2)=0$. So $X$ has one, respectively, three singularities,
according to $a_{21}^2-a_{12}a_{30}$ is negative, respectively
positive.  In both cases all the singular points of $X$ are
hyperbolic saddles and so, topologically, in a full neighborhood
of $0$ the implicit differential equation (\ref{eq:kgt}) is
equivalent to a Darbouxian umbilic point $D_1$ or to a Darbouxian
umbilic point of type $D_3$. See Fig. \ref{fig:d3p} and \cite{gs1,
gs5}.
\begin{figure}[htbp]
\par
\begin{center}
\includegraphics[angle=0, width=6.5cm]{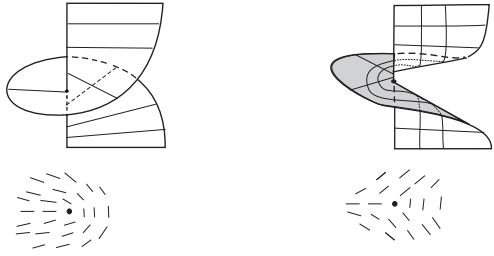}
\end{center}
\caption{ \label{fig:d3p} Resolution of a singular point by a Lie-Cartan line field }
\end{figure}

In fact,

$$\aligned DX(0,0,p)=& \left( \begin{matrix}  -2 a_{21}p-a_{30} &-2 a_{12}p-a_{21}  &0\\
-p(2 a_{21}p+a_{30}) &-p(2 a_{12}p+a_{21})  &0\\ A_{31} & A_{32}
&A_{33} \end{matrix} \right )\endaligned
$$
\noindent where, $A_{31}=p( (c +a_{22})p^2   +2 a_{31}p+a_{40})$,
  $A_{32}= p [ a_{13}p^2+( 2a_{22}-c)p+a_{31}]$ and $A_{33}=4 a_{21}p+a_{30}+3a_{12}p^2$.

The eigenvalues of $DX(0,0,p)$ are $\lambda_1(p)=-(a_{30}+
 3a_{21}p +2a_{12}p^2)$, $\lambda_2(p) =4 a_{21}p+a_{30}+3a_{12}p^2$ and $\lambda_3=0$.

Let $p_1$ and $p_2$ be  the roots of $ r(p)=  a_{30}+2
a_{21}p+a_{12} p^2 =0$.

Therefore,  $\lambda_1(p_i) = a_{21}p_i+a_{30}$ and
$\lambda_2(p_i) = -2(a_{21}p_i+a_{30})$.

As   $$r(-\frac{a_{30}}{a_{21}})=
\frac{a_{30}(a_{12}a_{30}-a_{21}^2)}{a_{21}^2}\ne 0,$$

\noindent it follows that $\lambda_1(p_i)\lambda_2(p_i) <0$ and
$\lambda_1(0)\lambda_2(0) = -a_{30}^2<0$. So the singularities of
$X$ are all hyperbolic saddles. If $  a_{21}^2 -a_{12}a_{30}<0$,
$X$ has only one singular point $(0,0,0)$. If $  a_{21}^2
-a_{12}a_{30}>0$, $X$ has   three  singular points $(0,0,0)$,
$(0,0,p_1)$ and $(0,0,p_2)$.

In the first case in a full neighborhood of $(0,0)$ the principal
foliations have the topological type of a $D_1$ Darbouxian umbilic
point. In the region $w>0$ the behavior is as shown in Fig.
\ref{fig:kgt}  (left). In the second case the principal foliations
have the topological type of a $D_3$ Darbouxian umbilic point and
so the behavior in the finite region $w>0$ is  as shown in Fig.
\ref{fig:kgt}  (right).
\end{proof}

\section{Principal Nets at Critical End Points}\label{sec:end2}

Let $p$ be a critical end point of the surface $A(\alpha)$,
$\alpha=\alpha^c\circ\mathbb P$. Without lost of generality assume
that the point $p$ is located at $(0,1,0,0)$  and that the surface
$\alpha^c =0$ is given by the graph of a function $w=h(u,v)$,
where $h$  vanishes together with its first partial derivatives at
$(0,0)$ and the $u$ and $v$   are the principal axes of the
quadratic part of its second order jet.

Through the central projection $\mathbb Q$, the coordinates $(u,v,w)$ can be thought to be orthonormal in the tangent space ${\mathbb T}
_p^3$ to $\mathbb S^3$ at $p$, with $w$ along
$\omega
=(0,0,0,1)$,   $u$ along $(1,0,0,0)$  and $v$ along $(0,1,0,0)$.

\begin{lemma} \label{lem:c2} Let $p$ be a critical end point of the surface $A(\alpha)$, $\alpha=\alpha^c\circ\mathbb P$.  Then  there is
  a mapping $\overline{\alpha}$ of the form
$$\overline{\alpha}(u,v)=(x(u,v), y(u,v), z(u,v))$$
 defined by

\begin{equation}\label{eq:cep} \aligned x(u,v)=&
\frac{u}{h(u,v)},\;\;\;  y(u,v)= \frac{v}{h(u,v)}, \;\;\; z(u,v)=
\frac{1}{h(u,v)}\endaligned \end{equation}
which parametrizes the surface $A(\alpha)$ near $p$.  The function $h$ is as follows.
\begin{enumerate}
\item[i)] If $p$ is  a     definite critical point of $ h$, then
\begin{equation}\label{eq:cepd}
\aligned
 h(u,v) =&  (a^2 u^2+b^2 v^2)+\frac 16(
a_{30}u^3+3a_{21}u^2v+3a_{12}uv^2+a_{03}v^3)\\
+&\frac 1{24}( a_{40}u^4+6a_{31}u^3v+4 a_{22}u^2v^2+ 6
a_{13}uv^3+a_{04}v^4)+ h.o.t. \endaligned \end{equation}

\item[ii)] If $p$ is  a   saddle  critical
point of   $h$, then
 \begin{equation}\label{eq:ceps}
\aligned
 h(u,v)=&  (-au +   v )v +\frac 16(
a_{30}u^3+3a_{21}u^2v+3a_{12}uv^2+a_{03}v^3)\\
+&\frac 1{24}( a_{40}u^4+6a_{31}u^3v+4 a_{22}u^2v^2+ 6
a_{13}uv^3+a_{04}v^4)+ h.o.t.
\endaligned
\end{equation}
\end{enumerate}

\end{lemma}

\begin{proof}
The map $x=u/w, \; y=v/w, \; z=1/w$ from  $\mathbb T
_y^3$  to $\mathbb T
_w^3$,  expresses the composition  ${\mathbb P  }^{-1}\circ {\mathbb Q    }$.

 Therefore the surface $A(\alpha)$, with  $\alpha = \alpha^c \circ (\mathbb P) ^{-1}$ can be parametrized with the functions
 $x, y, z$ as is stated in equation (\ref{eq:cep}).

 The function $h$ takes the form   given in equation (\ref{eq:cepd})  if it is
definite positive. If it is a non degenerate saddle, after a
rotation of principal axes, $h$ can be written in the form   given
in equation (\ref{eq:ceps}).
\end{proof}

\begin{lemma} \label{lem:d2} The differential equation (\ref{eq:11p}) in the chart
$\overline{\alpha}$ of Lemma \ref{lem:c2}, multiplied by
$h^4\sqrt{EG-F^2}$, extends to a full domain of the chart $(u,v)$
to one   given by
\begin{equation} \label{eq:lmns} \aligned
   L&dw^2+Mdudw+Ndu^2=0,  \\
L=& h^4\sqrt{EG-F^2}(Fg-Gf),\\
 M=& h^4\sqrt{EG-F^2}(Eg-Ge), \\
 N=&h^4\sqrt{EG-F^2}(Ef-Fe).
\endaligned
\end{equation}
\noindent where the coefficients are of   class $C^{k-2}$.
Here $(E,F,G) $ and $(e,f,g)$ are the coefficients of the first and second fundamental forms of the surface in the chart $\overline{\alpha}$.

\end{lemma}

\begin{proof}
 The first fundamental
form of the surface parametrized by $\overline{\alpha}$, equation
(\ref{eq:cep}),  in Lemma \ref{lem:c2} is given by:

$$\aligned
 E(u,v) =& \frac{ (h-uh_u)^2
 +( v^2+1)h_u^2  }{h^4}\\
F(u,v)=&\frac{-h(uh_u+v h_v)+(u^2+v^2+1)h_u h_v  }{h^4} \\
G(u,v)=& \frac{( h - vh_v)^2
 +(u^2 +1)h_v^2  }{h^4}
\endaligned
$$
 The coefficients of the second fundamental
form are given by :

$$\aligned
 e(u,v) =& -\frac{ h_{uu }}{h^4\sqrt{EG-F^2}}\\
f(u,v)=&-\frac{ h_{uv}  }{h^4\sqrt{EG-F^2}} \\
g(u,v)=&- \frac{h_{vv}}{h^4\sqrt{EG-F^2}}
\endaligned
$$
\noindent where $e=[\overline{\alpha}_{uu}, \overline{\alpha}_u, \overline{\alpha}_v]/|\overline{\alpha}_u\wedge \overline{\alpha}_v|, $  $
f=[\overline{\alpha}_{uv}, \overline{\alpha}_u, \overline{\alpha}_v]/|\overline{\alpha}_u\wedge \overline{\alpha}_v|  $ and $g=[\overline{\alpha}_{vv},
\overline{\alpha}_u, \overline{\alpha}_v]/|\overline{\alpha}_u\wedge \overline{\alpha}_v|. $

Therefore the differential equation of curvature lines, after
multiplication by $h^4|\overline{\alpha}_u\wedge \overline{\alpha}_v|$ is as stated.
\end{proof}

\subsection{Differential Equation of Principal Lines around a Definite  Critical End  Point} \label {sec:edfocal}

\begin{proposition} \label{prop:eqfoco} Suppose that $0$ is a critical point of $h$ given by equation (\ref{eq:cepd}),
 with
 $a>0, \; b>0$ (local   minimum).

In polar coordinates $u=br\cos\theta$, $v=ar\sin\theta$   the
differential equation (\ref{eq:lmns}) is given by
$Ldr^2+Mdrd\theta +Nd\theta^2=0$, where:

\begin{equation}\label{eq:lmnp} \aligned L=&l_0+l_1 r +h.o.t,\\
 \ M=&m_0+m_1r +h.o.t,\\
 N=&r^2 (\frac 12 n_0+\frac 16n_1r+\frac{1}{24}n_2r^2+h.o.t.)\endaligned
\end{equation}

\noindent with $ m_0=M(\theta,0)= -8a^7b^7\ne 0$ and the
coefficients $(l_0, l_1, m_1,$ $ n_0,$ $ n_1, n_2)$ are
trigonometric polynomials with coefficients depending on the
fourth order jet of $h$ at $(0,0)$, expressed in equations
(\ref{eq:l0}) to (\ref{eq:n2}).
\end{proposition}

\begin{proof}
Introducing polar coordinates $u=br\cos\theta$, $v=ar\sin\theta$
in the equation (\ref{eq:lmns}), where    $h$ is given by equation
(\ref{eq:cepd}), it follows that the differential equation of
curvature lines near the critical end point $0$,    is given by
$Ldr^2+Mdrd\theta +Nd\theta^2=0$, where: $m_0= M(\theta,0)=
-8a^7b^7, \;\; N(\theta ,0)= 0, $ and $\frac{\partial N}{\partial
r}(\theta, 0)=0$.

The Taylor expansions of $L$, $M$ and $N$ are as follows:
$$\aligned L=&l_0+l_1 r +h.o.t,\\
  M=&m_0+m_1r +h.o.t,\\
 N=&r^2 (n_0/2+n_1r/6+n_2r^2/24+h.o.t. )\endaligned
$$
  After a long calculation, corroborated by computer algebra, it follows that:

\begin{equation}\label{eq:l0} \aligned l_0=& =2 a^5 b^5 [a_{30} b^3\cos^2\theta \sin\theta +
a_{21}a
b^2(2\cos\theta -3  \cos^3  \theta)
\\+& a_{12}a^2b( \sin\theta - 3  \cos^2\theta  \sin\theta)
 +a_{03}
a^3(\cos^3\theta -\cos\theta)]
\endaligned
\end{equation}

\begin{equation}\label{eq:n0} \aligned n_0 =&4 b^5 a^5[(-3 a_{21}b^2 a+a_{03}a^3) \cos^3\theta
\\
+&(a_{30}b^3-3a_{12}ba^2)\sin\theta
\cos^2\theta \\+& (-a_{03}a^3+2a_{21}b^2a)\cos\theta +4ba^2
a_{12}\sin\theta] \endaligned
\end{equation}

\begin{equation}\label{eq:m1}\aligned
 m_1=&-4 b^5
a^5[(a_{30}b^3+a_{12}ba^2)\cos\theta+(a_{21}b^2a+a_{03}a^3)\sin\theta]
\endaligned \end{equation}

\begin{equation}\label{eq:n2}
\aligned  n_2 =&-2 a^3 b^3 [   ( a^3b^4 (24 a_{30}   a_{13}
 +6a_{03} a_{40} +108 a_{21} a_{22}+  72 a_{12}
 a_{31})\\
+& 6 a^7 a_{03}a_{04}
-18b^2 a^5 (  a_{21} a_{04}+4 a_{12}  a_{13}+ 2 a_{03}
a_22)
-6ab^6(3 a_{21}   a_{40}\\
+&4 a_{30} a_{31}))
\cos^7\theta
 + (a^4b^3(  6 a_{04} a_{30}+ 72 a_{21} a_{13}+ 108 a_{12}   a_{22}\\
+&24 a_{03}
  a_{31}) + 6 a_{40} b^7 a_{30}
- 18 b^5 a^2(    a_{40}
 a_{12} +4 a_{21}  a_{31} +2 a_{30}  a_{22})\\
 -&  6 a^6 b  (4 a_{03}  a_{13} +3 a_{12}
 a_{04})
) \sin\theta \cos^6\theta
 + ( 24 a^7 b^2 a_{03} + 72 b^6 a^3 a_{21}\\
 -&b^4 a^5(24  a_{03}+72  a_{21})
 +   a b^6 (50 a_{30} a_{31}+36 a_{21}  a_{40} )\\
 -& a^3b^4(  10 a_{03}  a_{40} +198 a_{21}
  a_{22} +46 a_{30}  a_{13}  +126 a_{12}   a_{31} )\\
+& b^2 a^3 ( 18a_{12}^2 a_{21} -6 a_{21}^2  a_{03} - 12 a_{30}   a_{12}   a_{03} )
+b^4 a (  3 a_{30}^2   a_{03}   \\
+& 6 a_{30} a_{12}   a_{21}
-9
a_{21}^3)
 + a^5 b^2 (54 a_{03} a_22 +114 a_{12} a_{13}+   30 a_{21}
a_{04})\\
+ &a^5(3 a_{03}^2 a_{21}-3 a_{12}^2  a_{03}) -8a_{03} a^7 a_{04}) \cos^5\theta\\
+& ( a^6 b( 22 a_{03}  a_{13} +18 a_{12}   a_{04}  ) + a^4b^5(72
a_{12} +24  a_{30}) -10 a_{40} b^7 a_{30}\\
-& a^4 b^3(  90 a_{21} a_{13} + 126 a_{12}
a_{22} +  8 a_{04}  a_{30} + 26 a_{03}  a_{31}   ) \\
+& a^2 b^5 ( 54 a_{30} a_22  +24 a_{40}  a_{12} +102 a_{21} a_{31}
 )
 + a^2 b^3( 12 a_{30}  a_{03} a_{21} \\
+&6 a_{12}^2
 a_{30}
 -18 a_{21}^2 a_{12})
-24 b^7 a^2 a_{30} +  3 b^5(a_{21}^2   a_{30}- a_{30}^2  a_{12} ) \\
-&a^4b ( 3 a_{03}^2   a_{30}+ 6 a_{21}  a_{12} a_{03} - 9 a_{12}^3) -72 a^6 b^3 a_{12}) \sin\theta \cos^4\theta  \\
+&  (-2 a_{03} a^7 a_{04}-24 a^7 b^2 a_{03} + 6 a^5(   a_{12}^2
a_{03} - a_{03}^2 a_{21} )- a b^6( 16 a_{21}  a_{40} \\
+& 24 a_{30}  a_{31})
 -a^5 b^2 (  12 a_{21}  a_{04}+18
a_{03} a_22  + 48 a_{12}  a_{13})\\
 +& a^5b^4 (24  a_{03}  + 96 a_{21}
)
 + a^3b^2( 15 a_{30}   a_{12}   a_{03}+
 12 a_{21}^2  a_{03}-27 a_{12}^2
  a_{21} )\\
+& a^3b^4(  4 a_{03} a_{40}  + 120 a_{21}  a_22     +70 a_{12}
  a_{31}  +26 a_{30}
a_{13}    )\\
+& ab^4 ( 6 a_{21}^3  -3 a_{30}^2 a_{03}- 3 a_{30}   a_{12} )
a_{21} -96 b^6 a^3 a_{21}
 ))
\cos^3\theta
+( 48 a^6 a_{21} b^3\\
+&(4 a_{13} a_{03}+2 a_{21} a_{04}) a^6 b + (9
a_{21} a _{03} a_{21}+3 a_{03}^2 a_{30}-12 a_{21}^3)
 a^4 b\\
-&( 48 a_{21}+24 a_{30}) a^4 b^5+ (6 a_{03} a_{31}+34 a_{13}
a_{21}+48 a_{21} a_{22}\\
+&2 a_{04} a_{30}) a^4 b^3
 +24 b^7 a^2 a_{30}
-(36 a_{31} a_{21}+18 a_{22} a_{30}+6 a_{40} a_{21}) a^2
 b^5\\
+&(9 a_{21}^2 a_{21}-9 a_{30} a_{03} a_{21}) a^2 b^3)
\sin\theta \cos^2\theta\\
 +&  (4 a_{03} a^7 a_{04}+(-3 a_{12}^2 a_{03}+3 a_{03}^2 a_{21}) a^5\\
 +&6 a_{13} b^2 a^5 a_{12}-36 a^5 b^4 a_{21}
 +3(3a_{12}^2 a_{21}-2 a_{21}^2 a_{03}-  a_{12} a_{03} a_{30}) b^2 a^3\\
 -&
 4( 6 a_{21} a_{22}+  a_{13} a_{30}+3 a_{12} a_{31}) b^4 a^3+36 b^6 a^3 a_{21}) \cos\theta\\
  + &( 12 b^5 a^4 a_{12}-12 a^6 b^3 a_{12}- ( 2 a_{12}a_{04}+ 2 a_{13}a_{03}) a^6 b\\
 -&( 4 a_{13}a_{21}-6 a_{12}a_{22})a^4b^3+(3a_{12}^3-3a_{12}a_{03}a_{21})a^4b ) \sin\theta].
\endaligned
\end{equation}

\begin{equation}\label{eq:l1}
\aligned l_1=& -  \frac 16 a^3 b^3[(18 b^5 a a_{21} a_{30}+18 b
a^5 a_{12} a_{03} \\
-& 6 b^3 a^3 (a_{30} a_{03}+9
 a_{21} a_{12})) \cos^6\theta \\
+&(3 a^6 a_{03}^2-9(2a_{21}a_{03}+3a_{12}^2) b^2
a^4\\
+&9(2a_{30}a_{12}+3a_{21}^2)b^4a^2-3b^6a_{30}^2) \sin\theta
\cos^5\theta  \\
+&(-15 b^5 a a_{21} a_{30}+9 b^3 a^3(9 a_{21} a_{12}+
a_{30} a_{03}) \\
-&39 b a^5 a_{12} a_{03}-32 b^3 a^5 a_{13}+32 b^5 a^3 a_{31})
\cos^4\theta \\
+&(48 b^4 a^4 a_{22}-8 b^6 a^2 a_{40}-6 b^4 a^2 (3a_{21}^2+2 a_{12}
a_{30}) -8 b^2 a^6
a_{04}\\
-& 6 a^6
a_{03}^2+12 b^2 a^4(2 a_{21} a_{03}+3 a_{12}^2  )) \sin\theta \cos^3\theta  \\
+&(-24 b^5 a^3 a_{31}+40 b^3 a^5 a_{13}+24 b a^5 a_{12} a_{03} \\
-&3
b^3 a^3(9 a_{21} a_{12}+ a_{30} a_{03})) \cos^2\theta  \\
+&(-12 b^6 a^4-3 b^2 a^4(3a_{12}^2 +2
 a_{21} a_{03})+8 b^2 a^6 a_{04}\\
 +&12 b^4 a^6+3 a^6 a_{03}^2-24 b^4 a^4 a_{22}) \sin\theta
\cos\theta \\
-&( 3 b a^5 a_{12} a_{03}+8 b^3 a^5 a_{13})]\endaligned
\end{equation}

\begin{equation}\label{eq:n1} \aligned
n_1=& a^3 b^3[(18 b^5 a a_{21} a_{30}-6  b^3 a^3(9 a_{21} a_{12}+
  a_{30} a_{03})
\\
+& 18 b a^5 a_{12} a_{03}) \cos^6\theta\\
+& (9 b^4 a^2(2 a_{12} a_{30}+3  a_{21}^2)-9 b^2 a^4(2
a_{21}a_{03}+3 a_{12}^2)\\
-& 3 b^6 a_{30}^2+3 a^6 a_{03}^2) \sin\theta \cos^5\theta \\
+& (16 b^3 a^5 a_{13}+9 b^3 a^3( a_{30} a_{03}+9   a_{21} a_{12})\\
-&16 b^5 a^3 a_{31}-39 b^5 a a_{21} a_{30}-15 b a^5 a_{12} a_{03}) \cos^4\theta \\
+& (-24 b^4 a^4 a_{22}+6 b^6 a_{30}^2-12 b^4 a^2(3 a_{21}^2+2   a_{12} a_{30})+4 b^2 a^6 a_{04}\\
+&4 b^6 a^2 a_{40}+6 b^2 a^4 (2a_{21} a_{03}+3 a_{12}^2)  ) \sin\theta \cos^3\theta \\
+&  (-20 b^3 a^5 a_{13}-6 b a^5 a_{12} a_{03}+12 b^5 a^3 a_{31}\\
+&18 b^5 a a_{21} a_{30}
- 3b^3 a^3(13   a_{21} a_{12}+ a_{30} a_{03})) \cos^2\theta \\
+&  (-12 b^6 a^4-3 a_{12}^2 b^2 a^4+12 b^4 a^4 a_{22}+12 b^4 a^6\\
+& 6b^4 a^2(2 a_{21}^2+  a_{12} a_{30})-3 a^6 a_{03}^2-4 b^2 a^6 a_{04}) \sin\theta \cos\theta\\
+& (4 b^3 a^5 a_{13}+6 b^3 a^3 a_{21} a_{12}+3 b a^5 a_{12}
a_{03})]. \endaligned
\end{equation}

\end{proof}

\subsection{Principal Nets around a Definite   Critical End  Points} \label {sec:focal}

\begin{proposition} \label{prop:foco} Suppose that $p$ is an end  critical point and consider the chart defined
in Lemma \ref{lem:c2} such that   $h$  is given by equation
(\ref{eq:cepd}), with   $a>0, \; b>0$ (local   minimum).
  Then the behavior of curvature lines near $p$ is the
following.
\begin{itemize}
\item[i)] One principal foliation is radial.

\item[ii)] The other principal foliation  surrounds   $p$ and  the associated    return
map $\Pi$ is such that $\Pi(0)=0,$  $\Pi^{\prime } (0)=1,$ $\Pi^{\prime\prime }
(0)=0, $ $\Pi^{\prime\prime\prime } (0)=0$ and
$\Pi^{\prime\prime\prime\prime} (0) =
\dfrac{\pi}{2^{10}a^5b^5}\Delta$, where
\end{itemize}

$$\aligned \Delta=&    12(  a_{30}  a_{21}+3  a_{03}  a_{30}-5 a_{12}  a_{21})  b^6  a^4\\
+& 12(5 a_{12}  a_{21}-   a_{12}  a_{03}
 -3   a_{03}  a_{30})  a^6  b^4\\
+&4(3 a_{04}  a_{30}  a_{21}+   a_{13}  a_{21}^2+
10  a_{31}  a_{03}  a_{21})a^4b^4 \\
-&4(10  a_{13}  a_{12}  a_{30}
+3  a_{40}  a_{12}  a_{03}+  a_{31}  a_{12}^2)  a^4  b^4\\
+&4( a_{13}a_{03}^2-   a_{04}  a_{12}  a_{03}   )  a^8
+4(  a_{40}
a_{30}  a_{21}- a_{31}
a_{30}^2)  b^8\\
+&3(  a_{30}^3  a_{03}+ 2 a_{30}
a_{21}^3-3  a_{30}^2  a_{21}  a_{12})  b^6\\
+&3(3  a_{12}  a_{03}^2  a_{21} -2  a_{12}^3  a_{03} -  a_{30}  a_{03}^3)  a^6\\
+&4[ a_{03}(2  a_{13}  a_{21}  -3  a_{04}   a_{30} -3 a_{31}  a_{03}+12 a_{22}  a_{12})\\
+&5 a_{04}
a_{12}  a_{21}-13
a_{13}  a_{12}^2]  a^6  b^2 \\
+&4[ a_{30}(+
  3  a_{13}  a_{30}-2  a_{31}  a_{12}   +3 a_{40}  a_{03}-12  a_{22}  a_{21})\\
 -&5
a_{40} a_{21}  a_{12}   +13  a_{31}  a_{21}^2]  a^2  b^6\\
  +&9(   a_{30}  a_{21}^2  a_{03}-2
a_{30}  a_{21}  a_{12}^2+   a_{12}  a_{21}^3)  a^2  b^4 \\
+&  9(- a_{12}^3 a_{21}-  a_{30} a_{03} a_{12}^2+2 a_{12} a_{21}^2
a_{03})  a^4  b^2\\
+& 12  b^2  a^8  a_{12}  a_{03}- 12  b^8 a^2
a_{30} a_{21}
\endaligned $$
\end{proposition}

\begin{figure}[htbp]

\par
\begin{center}
\includegraphics[angle=0, width=6.5cm]{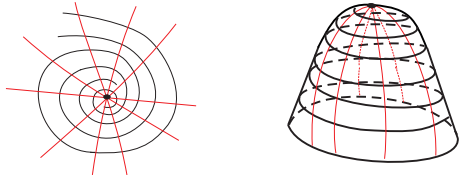}
\end{center}
\caption{ \label{fig:foco4} Curvature lines near a definite focal critical
  end point, $\Delta >0$. }
\end{figure}

\begin{proof}
Consider the implicit differential equation (\ref{eq:lmnp})
$$\aligned (l_0&+rl_1+h.o.t)dr^2+ (m_0+m_1r+h.o.t)drd\theta\\
+&r^2 (n_0/2+n_1r/6+n_2r^2/24+h.o.t. )d\theta^2=0.\endaligned $$

As $m_0=-8a^7b^7\ne 0$ this equation factors in the product
of two equations in the standard form, as follows.
\begin{equation}\aligned \label{eq:ret}
\frac{dr}{d\theta}=& -\frac 12 \frac{n_0}{m_0} r^2 + \frac 16
\frac{ 3 m_1-n_0+m_0 n_1 }{m_0^2 }r^3 \\
-&  \frac 1{24}\frac{
 (12 m_1^2 n_0+m_0^2 n_2+6 l_0 n_0^2-6 m_0 m_2 n_0-4 m_0 m_1 n_1)}{m_0^3}
 r^4\\
+&h.o.t.\\
 =& \frac 12 d_2(\theta) r^2+ \frac 16 d_3(\theta) r^3+ \frac 1{24} d_4(\theta)
 r^4 +h.o.t.
\endaligned
\end{equation}
\begin{equation}\label{eq:fr}
 \frac{d\theta}{dr}= -\frac {l_0}{m_0} +h.o.t.
\end{equation}

The   solutions of the nonsingular differential  equation
(\ref{eq:fr})
 defines the radial foliation.

Writing $r(\theta
,h):=h+q_1(\theta)h+q_2(\theta)h^2/2+q_3(\theta)h^3/6+q_4(\theta)
h^4/24+ h.o.t. $ as the solution of differential equation
(\ref{eq:ret}) it follows that:
\begin{equation}\label{eq:q}
\aligned   q_1^\prime (\theta)=&0 \\
 q_2^\prime (\theta)=& d_2(\theta)=-\frac{n_0}{m_0}  \\
q_3^\prime (\theta)=&  3 d_2(\theta) q_2(\theta)+ d_3(\theta) \\
q_4^\prime (\theta)=&
 3 d_2(\theta)q_2(\theta)^2+ 4 d_2(\theta)q_3(\theta)+ 6 d_3( \theta) q_2(\theta)+
 d_4(\theta)\endaligned
 \end{equation}
As $q_1(0)=0$ it follows that $q_1(\theta)=0$. Also $q_i(0)=0, \;
i=2,\; 3,\; 4$. So it follows that $q_2(\theta)=-\int_0^\theta
\frac{n_0}{m_0}d\theta$. From the expression of $n_0$, an odd
polynomial in the variables $c=\cos\theta$ and $s=\sin\theta$, it
follows that $q_2(2\pi)=0$ and therefore $\Pi^\prime(0)=1,
\;\;\Pi^{\prime\prime}(0)=1$.

Now, $q_3(\theta)=\int_0^\theta [q_2(\theta) q_2^\prime (\theta)
+d_3(\theta)]d\theta$.

Therefore $q_3(\theta)= \frac 12 q_2^2(\theta) +\int_0^\theta
d_3(\theta)d\theta.$

So,$$ \Pi^{\prime\prime\prime}(0)=
q_3(2\pi)=\int_0^{2\pi}d_3(\theta)d\theta.$$

A long calculation, confirmed by algebraic computation, shows that
$q_3(2\pi)=0$.

Integrating the last linear equation in (\ref{eq:q}),  it follows
that:

$$q_4(\theta)=3 q_2(\theta)^3+ 4 q_2(\theta)\int_0^\theta
d_3(\theta)d\theta+2\int_0^\theta
q_2(\theta)d_3(\theta)d\theta+\int_0^\theta d_4(\theta)d\theta.$$

Therefore,$$ \Pi^{\prime\prime\prime\prime }(0)= q_4(2\pi)=
2\int_0^{2\pi}
q_2(\theta)d_3(\theta)d\theta+\int_0^{2\pi}d_4(\theta)d\theta.$$

Integration of the right hand member, corroborated by algebraic computation, gives
  $ \Pi^{\prime\prime\prime\prime
}(0)=\dfrac{\pi}{2^{10}a^5b^5} \Delta. $
This ends the proof.
\end{proof}

\begin{remark} When  $\Delta \ne 0$ the foliation studied above spirals around $p$.
 The   point   is them  called a focal definite critical end point.

\end{remark}

\subsection{Principal Nets at Saddle Critical  End Points }\label{sec:saddle}

Let $p$ be a   saddle  critical
point of   $h$ as in equation
 (\ref{eq:ceps}) with the
finite region defined by $h(u,v)>0$.

 Then
the differential equation (\ref{eq:lmns}) is  given by

\begin{equation}\label{eq:cruz}\aligned
 L&dv^2+Mdudv+Ndu^2=0,\\
L(u,v)= & -a^3 u^2+2 a^2 u v -3 a a_{12}u v^2 + (2
a^2a_{12}-3 a a_{21}-2 a_{30})
  u^2v\\
  +& (a a_{30}+2 a^2 a_{21}) u^3 +(2 a_{12}+a a_{03})
 v^3+h.o.t.\\
 M(u,v)=& -2 a^2 v^2 +(4 a_{30}-a^2 a_{12})uv^2 +
 (a^2 a_{21}-2 a a_{30})  u^2v \\
 +& a^2 a_{30}u^3+(2 a a_{12}-a^2 a_{03}+4 a_{21})
 v^3+h.o.t.\\
 N(u,v)=& a v^2 [a^2-2 (a a_{21} + a_{30})
 u-2(aa_{12}u+a_{21})v]+h.o.t.
\endaligned
\end{equation}

\begin{figure}[htbp]

\par
\begin{center}
\includegraphics[angle=0, width=9.5cm]{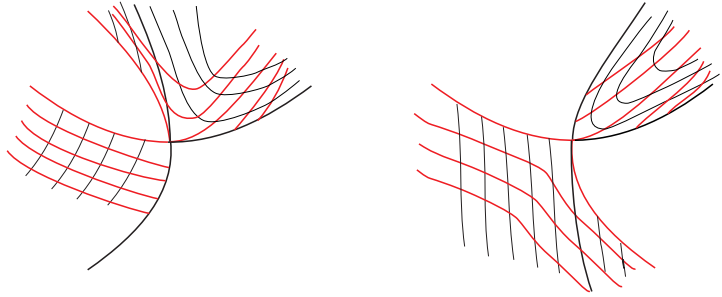}
\end{center}
\caption{ \label{fig:cruz} Curvature lines near saddle  critical
end point: even case, left,  and odd case, right. }
\end{figure}

\begin{proposition} \label{prop:saddle}
Suppose that $p$ is a saddle critical point of the surface
represented by $w=h(u,v)$ as in Lemma \ref{lem:c2}. Then the
behavior of the extended principal foliations
 in the  region ($h(u,v)\geq 0$), near  $p$, is the following.
\begin{itemize}
\item[i)] If $a a_{30}( a_{03} a^3  + 3
a a_{21} + 3 a^2 a_{12}  + a_{30}) >0$ then the curvatures of both  branches of   $h^{-1}(0)$ at
$p$ have the same sign and
 the behavior is as in Fig.
\ref{fig:cruz}, left -even case.

\item[ii)] If $a a_{30}( a_{03} a^3  + 3
a a_{21} + 3 a^2 a_{12}  + a_{30}) <0$ then the  curvatures of both  branches of   $h^{-1}(0)$ at
$p$ have opposite signs and the behavior is as in the Fig.
\ref{fig:cruz}, right - odd case.
\end{itemize}
\end{proposition}

\begin{proof} In order to analyze the behavior of the principal
lines near the branch of $h^{-1}(0)$ tangent  to $v=0$ consider
the projective blowing-up $ u=u, v=uw$.

 The
differential equation $Ldv^2+Mdudv+Ndu^2=$ defined by equation
(\ref{eq:cruz}) is,  after some simplification,  given by:

\begin{equation}\label{eq:cp}\aligned (-&\frac 14  a_{30}^2  u+a  a_{30}w +O(2)) du^2+(a  a_{30}u-2 a^2 w+O(2)) dw du \\
-&(a^2
u+O(2)) dw^2 =0.\endaligned
\end{equation}
To proceed consider the  resolution of the singularity $(0,0)$ of
equation (\ref{eq:cp}) by
 Lie-Cartan line field
 $X=(q{\mathcal G}_q, {\mathcal G}_q,-( q{\mathcal G}_u+{\mathcal G}_w)),\;\; q=\frac{du}{dw}. $
Here $\mathcal G$ is
 $${\mathcal G}=(-\frac 14  a_{30}^2  u+a  a_{30}w +O(2)) q^2+(a  a_{30}u-2 a^2 w+O(2)) q-(a^2
u+O(2)).   $$

The singularities  of $X$, contained in axis $q$ (projective
line), are the solutions of  the equation $q (2
a-a_{30}q)(6a-a_{30}q)=0$.

Also,
$$\aligned DX(0,0,q)=& \left( \begin{matrix}  \frac 12({ q a_{30}(2a-a_{30}q)}) &-2 q a (a-a_{30}q ) &0\\
\frac 12({   a_{30}(2a-a_{30}q)})&-2  a (a-a_{30}q ) &0\\ 0&0&
A_{33}\end{matrix} \right )\endaligned
$$
\noindent where $A_{33}= 3 a^2-4 a a_{30}q+\dfrac 34 a_{30}^2 q^2.$

 The eigenvalues of $DX(0,0,q)$
are $ \lambda_1(q) =-2 a^2+3 a a_{30}q-\dfrac 12 a_{30}^2 q^2$,
$\lambda_2(q) = 3 a^2 - 4 a a_{30}q+ \dfrac 34 a_{30}^2 q^2$ and
$\lambda_3(q) =0$.

Therefore the non zero eigenvalues of $DX(0)$ are $-2a^2$ and
$3a^2$. At $q_1=\frac{2a}{a_{30}}$ the eigenvalues of
$DX(0,0,q_1)$ are $2a^2$ and $-2a^2$. Finally at
$q_2=\frac{6a}{a_{30}}$ the eigenvalues of $DX(0,0,q_2)$ are
$6a^2$ and $-2a^2$.

As a conclusion of this analysis we assert that the   net of
integral curves of equation (\ref{eq:cp}) near $(0,0)$ is the same
as one of the generic singularities of quadratic differential
equations, well known as the Darbouxian $D_3$ or a tripod,
\cite{gs1, gs5}. See Fig. (\ref{fig:d3p}).

  Now observe  that the curvature at $0$ of
 the branch of $h^{-1}(0)$ tangent  to $v=0$ is precisely $k_1=
 \dfrac{a_{30}}{a}$ and that  $h^{-1}(0)\setminus
 \{0\}$ is solution of equation (\ref{eq:cruz}). So, after the
 blowing down, only one branch of the invariant curve
 $v=\dfrac{a_{30}}{6a}u^2+O(3)$ is  contained in the finite region
 $\{(u,v): h(u,v)>0\}$.

 Analogously,  the analysis of the behavior of the principal lines near  the branch of $h^{-1}(0)$ tangent  to $v=au$
 can be reduced to the above case. To see this   perform   a rotation of angle
 $ \tan\theta= a$ and
 take new orthogonal coordinates $\bar u$ and $\bar v$ such that
 the axis $\bar u$ coincides with the line $v=au$.

 The curvature at $0$ of
 the branch of $h^{-1}(0)$ tangent
  to $v=au $ is
$$k_2=
- \dfrac{a_{03} a^3  + 3 a a_{21} + 3 a^2  a_{12}  + a_{30}}{3a}.$$

Performing the blowing-up $v=v, \;\; u=sv$ in the  differential
equation  (\ref{eq:cruz}) we conclude that it factors in two
transversal regular foliations.

Gluing the phase portraits studied so far and doing their blowing down, the net
  explained  below is obtained.

 The finite region $(h(u,v)>0)$ is formed by two sectorial regions
$R_1$, with  $\partial R_1= C_1\cup C_2$ and     $R_2$ with $\partial
R_2= L_1 \cup L_2.  $  The two regular branches of $h^{-1}(0)$ are given by $C_1\cup L_1$ and $C_2\cup L_2$.

If $k_1k_2<0  $ -- odd case --  then  one region, say $R_1$, is convex and
$\partial R_1$ is invariant for one extended principal foliation and
$\partial R_2$ is invariant for the other one. In each
region, each foliation has an invariant separatrix tangent to the
branches of $h^{-1}(0)$. See Fig. \ref{fig:cruz}, right.

If $k_1k_2> 0  $ -- even case -- then in a   region, say $R_1$,  the  extended principal
foliations  are equivalent  to a trivial ones, i.e., to $dudv=0$,
with  $C_1$ being a leaf of one foliation and $C_2$ a leaf of  the
other one. In the   region $ R_2$ each extended principal foliation has a
 hyperbolic sector, with separatrices tangent to the
branches of $h^{-1}(0)$ as shown in Fig. \ref{fig:cruz}, left.  Here
$C_2\cup L_1$ are leaves  of one principal foliation and $C_1\cap
L_2$ are leaves  of the other one.
\end{proof}

\section{Concluding Comments and Related Problems} \label{sec:fin}

We have studied here the simplest  patterns  of principal
curvature lines at end points,   as the supporting  smooth surfaces  tend
to infinity in $\mathbb R ^3$,  following  the paradigm established in \cite{gasalg} to describe
the  structurally stable patterns   for principal curvature
lines escaping to infinity on algebraic surfaces.

We have recovered here --see  Proposition  \ref{prop:he} -- the
main results of the structurally  stable inflexion ends
established in \cite{gasalg} for algebraic surfaces: namely the
hyperbolic and  elliptic cases.

In the present context a surface  $A(\alpha ^c)$ with $\alpha
^c\in {\mathcal A}_c^k$ is said to be {\it structurally stable} at
a singular end point $p$ if the $C^s,  $ topology with $s\leq k$
if the following holds. For any sequence of functions $\alpha_n^c
\in {\mathcal A}_c^k$ converging to $\alpha ^c$ in  the $C^s$
topology, there is a sequence $p_n$ of end points of
$A(\alpha_n^c)$ converging to $p$ such that the extended principal
nets of $\alpha_n=\alpha_n^c\circ \mathbb P$, at these points, are
topologically equivalent to extended principal net of
$\alpha=\alpha ^c\circ \mathbb P$, at $p$.

Recall (see \cite{gasalg}) that two nets $N_i \, i= 1, 2\,$  at
singular points $p_i \, i= 1, 2\,$ are topologically equivalent
provided there is a homeomorphism of   a neighborhood of  $p_1 $
 to  a neighborhood of  $p_2 $ mapping the respective points
 and
 leaves of the respective foliations to each
 other.

 The analysis in Proposition \ref{prop:he} makes
clear that
 the hyperbolic and elliptic inflexion end points are also structurally
  stable in the  $C^3$ topology for
 defining  $\alpha ^c$ functions in the space ${\mathcal A}_c^k$, $k\geq 4$.

 We have  studied also six new cases --see Propositions \ref{prop:kge} to \ref{prop:saddle} --  which represent the
simplest patterns  where the structural  stability conditions fail.

The lower  smoothness  class $C^k$ for the validity of the
analysis in the proofs of these propositions is as follows.  In
Propositions \ref{prop:kge} and \ref{prop:saddle}  we must assume
$k\geq 4$.  In Proposition \ref{prop:foco}, clearly $k\geq 5$ must
hold.

 In each of these cases
  it is not difficult to describe partial aspects of  possible
   topological changes --bifurcation
phenomena-- under small perturbations of the defining functions
$\alpha ^c$.

However it involves considerably technical work to provide  the
full analysis of bifurcation diagrams of singular end points and
their  global effects in the principal nets.

 We recall here  that the
study of the bifurcations  of principal          nets  away from
end points, i.e.,  in  compact regions  was  carried out in
\cite{G-S4}, focusing the  umbilic singular points. There was also
established the connection between umbilic codimension one
singularities and their counterparts in critical points of
functions and the singularities of vector fields, following the
paradigm of first order structural stability in the sense of
Andronov and Leontovich \cite{al}, generalized and extended by
Sotomayor \cite{sihes}. Grosso modo this paradigm aims to
characterize the structurally stable singularities under small
perturbations inside the space of non- structurally stable ones.

To advance  an idea of the bifurcations at end points, below we
will suggest  pictorially the local bifurcation diagrams in the
three regular cases studied so far.

\begin{figure}[htbp]
\par
\begin{center}
\includegraphics[angle=0, width=10.5cm]{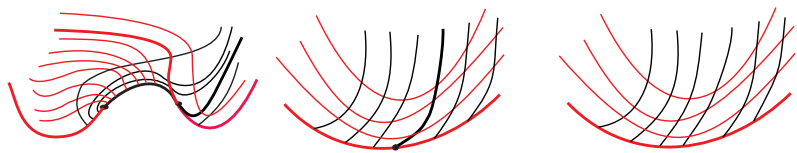}
\end{center}
\caption{ \label{fig:kgeu} Bifurcation Diagram of Curvature lines
near  regular end
 points: elimination of hyperbolic and elliptic  inflexion points }
\end{figure}

\begin{figure}[htbp]

\par
\begin{center}
\includegraphics[angle=0, width=10.5cm]{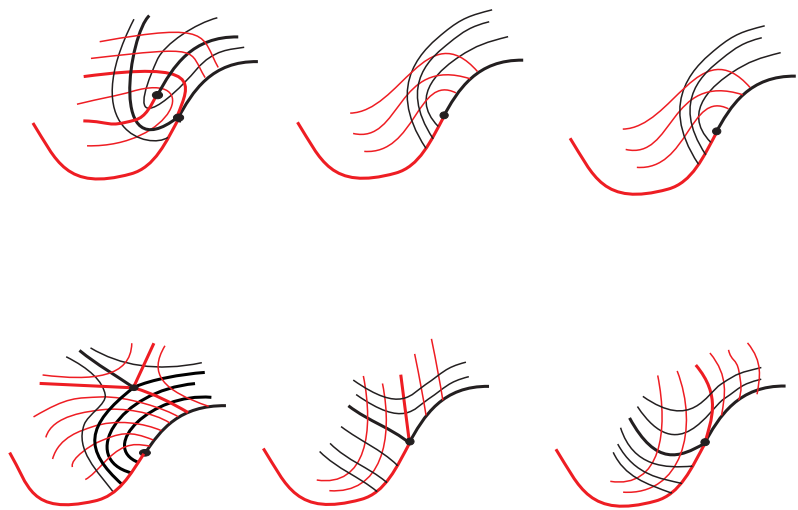}
\end{center}
\caption{ \label{fig:kgtu} Bifurcation Diagram of  Curvature lines
near  umbilic-inflexion   end
 points. Upper row: $D_1$ umbilic - hyperbolic inflexion.
Lower row: $D_3$ umbilic - elliptic inflexion.   }
\end{figure}

 The  description  of the bifurcations  in the  critical cases, however,
 is much more intricate  and will not be discussed
here.

The full analysis of the non-compact bifurcations as well as their
connection  with  first order structural stability  will be
postponed to a future paper.

Concerning the study of end points, see also   \cite{gs3}, where
Gutierrez  and Sotomayor studied  the behavior of principal nets
on constant  mean curvature surfaces, with special analysis of
their periodic leaves,  umbilic and  end points. However, the
patterns of behavior for this class of surfaces is non-generic in
the sense of the present work.

We conclude proposing the following problem.

\begin{problem}
Concerning the case of the focal critical end point,  we propose
to the reader to  provide a  conceptual analysis and a proof of
Proposition \ref{prop:foco},  avoiding long  calculations  and the
use of Computer Algebra.
\end{problem}

\end{document}